\documentclass[12pt, a4]{article}

\usepackage{amsmath,amssymb, amsfonts}

\title{\bf How to Make a New Logic\footnote{This is the first version of the paper. Some results 
will be added in the future version.}}
\author{Takao Inou\'{e}}
\date{August 13, 2021}

\begin{document}

\maketitle

\begin{abstract}
We discuss about how to make a new logic with considering a slogan: \it take the inversion of 
what you know. \rm  Suppose that we have a Gentzen-style logical system $S$. The operation to 
make a new logic from $S$ is the following: \it for all the axioms and rules of $S$, change the direction 
of all the arrows occurred in the sequents of the axiom or the rule to the opposite side. \rm 
We call this operation \it Stahlization. \rm We consider certain logics in this respect. 
\end{abstract}

\smallskip


\normalsize 

\section{The slogan for how to make a new logic}

   When we look at past revolutional works especially for philosophy, art, mathematics and 
science, which gave us totally different perspectives over individual field, such as 
Kant's Copernican revolution, J. S. Bach's 'Die Kunst der Fuge' (BWV 1080)\footnote{I mean pairs of 
Contrapunctus recus and inversus contrained in the immortal masterpies, where one sees 
beutiful symmetries and musical effects.}, 
homological algebra\footnote{In this branch of algebra, homology and cohomology play an important role.}, 
H. M. Friedman's reverse mathematics\footnote{H. M. Friedman's reverse mathematics can be summerized 
as the following themes of Friedman \cite{Friedman1975}: (i) when the theorem is from the right axioms, the axiom can be 
proved from the theorem; (ii) much more is needed to explicitly define a hard-to define set of integers that 
to merely prove their existance. For this current research area, see Simpson \cite{Simpson1987, Simpson2009}, 
Tanaka \cite{Tanaka1990} and so on.}, the Copernican system, Einstein's 
special theory of relativity and so on, we will be soon aware of some similarlity comprised in 
them\footnote{Rather technical examples are e.g. Lebesgue integration, the Boyer-Moore string-search algorithm, 
Inverse kinematics and so on.}. 
The similarity which we mean can be expressed as the following slogan:

\smallskip 

\centerline{\it Take the inversion of what you know. \rm } 

\smallskip

\noindent This is a very simple principle. However, the simple operation has already provided 
us with incredibly great novelties in spite of its simplicity\footnote{In this respect, B. C. Fraassen's 
book \cite{Fraassen} is interesting.}

There actually are many ways to make a new logic. We shall roughly classify the possible ways 
to two categories.

    \it The first category: \rm a modification of a given logic, in a broad sense, such as the symmetrization of 
the logic, some addition of new axioms or rules to the logic, taking a meaningful subsystem 
of the logic and so forth.

    \it The second category: \rm a creation of a totally new logic from a own philosophical point of view, 
or as a systematic description of some logical struture appeared in particular phenomena which 
we see. 

\noindent These two categories are not exclusive. For example, intuitionistic logic whose creating 
process would belong to the second category is nothing but a subsystem of classical logic as the 
result of the axiomatization of what L. E. J. Brouwer intended.

\section{Stahlization}

We think that we would have been enough prefaced with the above paragraphs. The aim of this 
essay is to make the readers realize that a simple operation to be mentioned below, satisfying the 
slogan introduced in the first paragraph and the first category just mentioned, is a very promising 
way to create a new logic. Suppose that we have a Gentzen-style logical system $S$. The operation to 
make a new logic from $S$ is the following: \it for all the axioms and rules of $S$, change the direction 
of all the arrows occurred in the sequents of the axiom or the rule to the opposite side. \rm We shall 
call this operation Stahlization after the name of a philosopher Gerold Stahl, mimicking 
'G\"{o}delization' in the terminology of recursion theory\footnote{We may call the operation \it dualization \rm 
instead of Sthalization. However, this terminology is not much fun.}.

We can find the example of Stahlization in two papers, i.e. Goodman \cite{Goodman1981} and 
Inou\'{e} \cite{Inoue1989}. Though \cite{Goodman1981} is ornamented by many considerations about his subject, it is, in a 
word, a study of the Stahliaztion of intuitionistic propositional calculus in Gentzen-style (for 
that, see Kleene \cite{Kleene1952}). He called the Stahlization of the logic \it anti-intuitionistic propositional 
calculus. \rm According to him, the choice of a logic in order to express his motivation and 
intention was due to F. W. Lawvere, one of great figures in the area of category theory\footnote{Better 
to say, in the area of the foundation of mathematics, we do not mean however that category theory 
is the foundation of mathematics.}. Goodman's motivation lies on a study in logic subject to 
vagueness or indeterminateness in the scope of \it  his \rm way and interest. His attempt is 
however still in the stage of try and error. We shall mention 
Goodman's contribution again in the sequel. Now, we shall proceed to another exmple of 
Stahlization. Stahl \cite{Stahl1958} proposed some propositional calculi to derive only the negation of a 
thesis in classical propositional logic. Such a system is called \it Stahl's opposite system\footnote{Some 
systems equivalent to Stahl's opposite system were independently studied by Morgan \cite{Morgan1973} and 
Bunder \cite{Bunder1989}.} \rm  and we 
call a thesis in such a system a \it contradiction. \rm The real novelty of Stahl's formulation from the 
standpoint of axiomatic rejection consists in the axiomatization of all the contradictions free from 
the notion of provability\footnote{Skura \cite{Skura1990} obtained some general result in theis respect. 
See also Inou\'{e} \cite{Inoue1989a}.}. In our use of words, Inou\'{e} \cite{Inoue1989} observed that the idea of Stahl is 
nothing but the Stahlization of classical propositional logic in Gentzen-style. His observation is 
based on the inversion principle after P. Lorenzen with respect to negation in Gentzen's sequent 
calculus LK\footnote{The inversion principle with respect to negation reads that if 
$\Gamma_1, \Gamma_2 \rightarrow \Theta_1, \sim A, \Theta_2$ holds in LK, then so does 
$\Gamma_1, A, \Gamma_2 \rightarrow \Theta_1, \Theta_2$.}. He could extend Stahl's logic 
to the first-order predicate level on the basis of the 
Stahlization (see \cite{Inoue1989} ). It should be mentioned that I, the author of \cite{Inoue1989} , also was stimulated by 
category theory to get the idea\footnote{Goodman's system is a consistent in the sense that 
not every sentence is derivable in the system.}.

    We remark that Goodman's anti-intuitionistic propositional calculus and my systems $SP$ and $SC$, 
in the notion of \cite{Inoue1989}, which are thus the Stahlization of classical propositional and predicate 
logics, since the former is a consistent system\footnote{I was also inspirated by studying the method 
of axiomatic rejection. For axiomatic rejection, see for example Inou\'{e} \cite{Inoue1989a} 
and Inou\'{e}, Ishimoto and Kobayashi \cite{Inoue2021}. Also refer to Carnielli and Pulcini \cite{Carnielli2017}, and 
Goranko, Pulcini and Skura \cite{Goranko2020}.} in which not all but some contradictions are 
derived, and the later is the extreme case of paraconsistent logic, namely the logic in which 
all theses are inconsitent. 

\section{Paraconsitent logic, Vasiliev's logic and beyond}

It is not a wonder that some of the current interests in philosophy and 
logic are in something vague or inconsistent such as vagueness, fuzzy logic, paraconsistent logic\footnote{As a recent literature, 
see Ba\c{s}kent and Ferguson \cite{Baskent2019}.} and so on. Also we know a tradition which goes back to Heraclitus, who insisted 
that the world is not consistent\footnote{See for example Rescher and Brandom \cite[pp. 1-2]{Rescher1980}.}. 
We also have to consider the influence of computer science to 
philosophy and logic. Computer science needs much more knowledge of inconsistent reasoning such 
as that of paraconsistent logic. Bodies of data stored in a computer memory may be not consistent, 
and if the stock of data is very large, it is not practical (thus not economical) to check the consistency 
of it at all\footnote{See da Costa and Marconi \cite[p. 24]{Costa1989}.}.

In connection with this essay, it is very important that G. Priest pointed out the similarlity between 
the interpretation of Vasiliev's neo-syllogistic proposed by Smirnov \cite{Smirnov1989} and that of anti-intuitionistic 
propositional calculus in Goodman \cite{Goodman1981} in his review Priest \cite{Priest1989} of \cite{Smirnov1989}. 
N. A. Vasiliev considered in the early 
years of this century atomic propositions of the form: $S$ is and is not $P$\footnote{For Vasiliev's 
works, see especially Smirnov \cite{Smirnov1986, Smirnov1989} and Arruda \cite{Arruda1977, Arruda1980, Arruda1984}.}. 
In the interpretion presented in Smirnov \cite{Smirnov1989}, such statements are true if all the points in the set $S$ 
are on the topological boundary of the set $P$. From this fact, I am tempted to say that there is a meaningful syllogistic and 
Vasiliev's logic is the Stahlization of the logic. To find such an original logic is kind of inverse problem. From the 
definition of Stahlization, if we have a Genzen-style formulation of Vasiliev's neo-syllogistic, 
we immediately obtain the original logic by taking the Stahlization of the neo-syllogistic. We will 
obtain much more about the neo-syllogistic by the simple operation.

If we are given a logic, the Stahlization of it will give us new perspective over the original logic, if 
the Gentzen-style formulation of the logic is possible. This means the creation of a new logic in 
other words. This is what I wanted to write in this essay, that is, one of my paradigm over logic.

\bigskip
$$ $$
$$ $$

\noindent \bf Acknowledgments. \rm I would like to take this opportunity to express my gratitude and deep condolences to 
the following three scholars who have passed away.

\medskip 

Logician, Marcin Mostowski (1955-2017),

Computer scientist, Grzegorz Bancerek (1966-2017),

Mathematician/Logician, Stanis\l{}aw \'{S}wierczkowski (1932-2015)

\bigskip

\noindent I would also like to thank an anonymous referee for suggestions on this paper.


\noindent Meiji Pharmaceutical University

\noindent Department of Medical Molecular Informatics

\noindent Tokyo, Japan 

\bigskip

\noindent Hosei University

\noindent Graduate School of Science and Engineering

\noindent Tokyo, Japan 

\bigskip

\noindent Hosei University

\noindent Faculty of Science and Engineering

\noindent Department of Applied Informatics

\noindent Tokyo, Japan

\medskip

\noindent ta-inoue@my-pharm.ac.jp

\noindent takao.inoue.22@hosei.ac.jp

\noindent takaoapple@gmail.com

\end{document}